\documentclass[a4 paper, 11pt]{article}
\usepackage[active]{srcltx} 
\usepackage[english]{babel}
\usepackage[T2A]{fontenc}
\usepackage[cp1251]{inputenc}
\usepackage[centertags]{amsmath}
\usepackage{amsfonts}
\usepackage{amssymb}
\usepackage{amsthm}
\usepackage{amsmath}
\usepackage{amscd}
\usepackage{latexsym}
 \newtheorem
   {theo}{Theorem}
\newtheorem {cit}{Theorem}
\title{Planar Lebesgue measure of exceptional set in approximation of subharmonic functions}

\author{Markiyan Girnyk (=Hirnyk) \\ Lviv Academy of Commerce \\Ukraine, 79005 Lviv, 10, Tuhan-Baranovskyi St.
\\
hirnyk@lac.lviv.ua }
\date{}
\begin{document}

\maketitle \begin{abstract}
We consider the
pointwise approximation of a subharmonic function by the logarithm
of the modulus of an entire one up to a bounded quantity.
 In the case of finite order an estimate from
below of the planar Lebesgue measure of an exceptional set in such approximation is
obtained.


\end{abstract}\large

Results on the approximation of a subharmonic function by the logarithm of the modulus of an entire function have numerous applications in complex analysis and potential theory
 (see, for example
 ,~\cite{A}-~\cite{E}). Such pointwise approximation is possible only outside an exceptional set, and for this reason the principal question concerning its minimal size arises. In this article we prove that the planar Lebesgue measure of an exceptional set in the approximation of the subharmonic function $|z|^\rho$ by the logarithm of the modulus of an entire function, having at most order $\rho$ and normal type, cannot be arbitrary small in some sense.

We use principal results and standard notations of potential theory ~\cite{F} and theory of value distribution ~\cite{K}. Let us recall some of them. We denote
 $D(a, r):=\{z:|z-a|<r\}
 ,\,\,C(a, r):=\{z:|z-a|\le r\}
 ,\,\,S(a, r):=\{z:|z-a|=r\},\,\,A(t, T]:=\{z:t< |z|\le T\},\,m_d$--  the Lebesgue measure on ${\mathbb{R}^d}$,
 letters $C$ with indices stand for positive constants, in parentheses we indicate dependance on parameters. As usually, $a^+= \max (a,0),\,a^- = \max (-a,0)$.
 Let $u$ be a subharmonic function, then
     $\mu_u$ is its Riesz measure,  $B(r,u):=\max\{u(z):z\in C(0,r)\}\,$ is
    the maximum,
    $n(a,r,u):= \mu(C(a,r)),\,n(r):=n(0,r,u)$ are the counting functions of the Riesz measure, $h(z,u,\mathcal{D})$ is the minimal harmonic majorant of the function
     $u$ in a domain $\mathcal{D}$
    (we sometimes omit it in notations),
    $T(r,u):=\frac 1{2\pi}\int_0^{2\pi}u^+(r e^{i \varphi})\,d \varphi$ is the Nevanlinna characteristic
    of $u$. It should be noted that the Nevanlinna characteristic of a meromorphic function $f$ is also denoted by $T(r,f)$.
    Since it is clear from the context which characteristic is used, this does not make the reader's difficulties.
    We denote by  $\Lambda$ the class of nondecreasing slowly changing functions
 $ \lambda
:[1,\infty) \mapsto [1,\infty)\,$(in particular, $\lambda (2r)\sim
\lambda (r)$ if $r\rightarrow \infty$).

   The notation $a\asymp b$ means that $|a|\le \mathrm{const}\cdot|b|$ and $|b|\le
    \mathrm{const}\cdot |a|$. Let  $\Theta\supset\Lambda$ be the class of nondecreasing functions $ \lambda
:[1,\infty) \mapsto [1,\infty)$ having the property
    : $\lambda(2R)\asymp \lambda (R)$. This implies the finitness of order
    $$ \tau=\limsup\limits_{R\to \infty} \frac {\log \lambda (R)}{\log R}$$
   of the  function $ \lambda $.
   The content of this work is closely associated with a theorem by I. Chizhikov  ~\cite{G}, which
   strengthens and specifies  a result by Yu. Lubarskii and Eu. Malinnikova
   ~\cite{J}, and with a theorem by R. Yulmukhametov ~\cite{L}.
     For the reader's convenience we cite these ones (in somewhat modificated, but equivalent formulations). \begin{cit}
     Let $u$ be a subharmonic function with the Riesz measure $\mu_u$. If for a function
    $\lambda \in \Lambda$ there exists a number $R_0$ such that for every number $R>R_0$ the condition
    \begin{equation} \mu_u(A(R,\,R\lambda(R)])>1\,
     \end{equation} holds, then there exists
     an entire function $f$ and a constant $C_1$ satisfying \begin{equation}
     \int_{A(R, 2R]}|u(z)-\log|f(z)||\,dm_2(z)<C_1R^2\log \lambda
     (R).
     \end{equation}
     Moreover, for every real number $\varepsilon > 0$ there exists a constant
     $C_2(\varepsilon)$ and a set $E(\varepsilon)$ such that \begin{equation}
|u(z)-\log|f(z)||< C_2(\varepsilon)\log \lambda (|z|),\,z\notin
E(\varepsilon),
     \end{equation} and \begin{equation}
     m_2(E(\varepsilon)\cap A(R,2R])/R^2<\varepsilon.
     \end{equation}\end{cit}
\begin{cit} Let $u$ be a subharmonic function of finite order $\rho$ and
a number $\alpha > \rho.$ Then there exists an entire function  $f$,
a constant $C_3=a_0+a_1\alpha,\,a_1 > 1,$, depending only on
$\alpha$ , and an exceptional set  $E$ depending on the functions
 $u,\,f$ and the number $\alpha$  such that
\begin{equation}    \left| u(z)-\log |f(z)| \right| \le C_3(\alpha) \log|z|
    ,\, z\notin E,\end{equation}
 where $E\subset\cup_j D(z_j,r_j) $ and
\begin{equation} \sum_{R<|z_j|\le 2R }r_j=o(R^{\rho -\alpha}),\quad R \to \infty.
     \end{equation}
\end{cit}
Let us formulate the result of our work. \begin{theo}
For all the real numbers
$\varepsilon>0, \rho>0$, for every an entire function
 $f$ satisfying the condition\begin{equation}
B(r,\log |f|)< C_4 r^\rho,
     \end{equation} for every function $\lambda \in \Theta$ of order $\tau$, and for any measurable set
      $E$,
      the condition
      \begin{equation} \left||z|^\rho-\log|f(z)|\right    |< C_5\log \lambda(|z|),\,z\notin E, \end{equation}
      implies the existence of a constant
       $C_6(\varepsilon)$
    such that\begin{equation} m_2(E\cap A(R,
     2R])>C_6(\varepsilon)R^{\chi+\rho},\,R>1.
     \end{equation}where $\chi=\min(2-2\rho-\varepsilon, 2-\rho-4C_5\tau-\varepsilon)
    .$ \end{theo}

Because the formulation of Theorem 1 is long, we formulate its statement in the important case of a bounded function
 $\lambda$ (its order
 $\tau=0$ and $\chi=2-2\rho-\varepsilon):$
$$m_2(E\cap A(R,
     2R])>C_6(\varepsilon)R^{2-\rho-\varepsilon},\, R>1.
$$

     Let us comment the content of Theorem 1. It states that the number
       $\varepsilon$ on the right-hand side of (4) cannot be replaced by an arbitrary function
      $\varepsilon(R)\rightarrow 0, \,R\to \infty$. Next, the condition (7) seems to be natural because it holds in Theorem 1 for subhamonic functions of finite order
     $\rho$ and normal type 1, following from
       (2). Indeed, let us consider the inequalities
      $$
      2\pi
      T(r,f)\le\int\limits_0^{2\pi}|\log|f(re^{i\varphi})||\,d\varphi+O(1)
      \le\int\limits_0^{2\pi}(|\log|f(re^{i\varphi})|-u(re^{i\varphi})|+u(re^{i\varphi}))\,d\varphi+O(1),
      $$which imply
      $$
      2\pi R^2T(R,f)<2\pi\int\limits_R^{2R}T(r,f)r\,dr\le
      \int_{A(R,2R]}|\log|f(z)|-u(z)|\,dm_2(z)+$$
      $$+\int_{A(R,2R]}u(z)\,dm_2(z)+O(R^2)\le C_1 R^2\log
      \lambda(R)+O(R^{\rho+2})+O(R^2),\,$$ $$R\to
      \infty,
      $$
and the growths of the functions $B(r,\log|f|)$ and  $T(r,f)$ are equal. We also note the possibility of a good approximation of a subharmonic function with finite order by the logarithm of the modulus of some entire function with infinite order. We present a simple example.
 Let
 $$\left||u(z)-\log|f(z)|\right |< C_5\log \lambda(|z|),\,z\notin E,
 $$
where $u$ is a subharmonic function of finite order, $f$ is an entire function, then

  $$\left||u(z)-\log|f(z)\cdot V(z)|\right |< C_5\log \lambda(|z|)+2,\,z\notin E\bigcup S:=\{z=x+iy:x\ge 1,
   |y|\le\pi \},
 $$
where (see~\cite[С. 256-258]{K} ) the function $$V(z)=\left\{
\begin{array}{l}\exp(\exp z)+\Psi_1(z)/z,\,z \in S;\\
\Psi_2(z)/z,\,z \not\in
 S;\\ \end{array} \right.
 $$ а $|\Psi_j(z)|<2,\,|z|>r_0>1,\,j=1,2.$
  We also note   that in
~\cite{L} and ~\cite{H}  the estimates from below of the sum of radii for any disk covering of the exceptional set are obtained, but no estimate from below of the planar measure  follows from those estimates.

 We cite two above mentioned results.
 \begin{cit}Let a number $\varepsilon>0$ and an entire function $f$ satisfy the inequality
 $$||z|-\log|f(z)||=o(\log |z|),\,\,E\not\ni z \to \infty,
 $$
  where $E\subset \bigcup_j\{z:|z-z_j|<r_j\},$  and radii $r_j$ are uniformly bounded. Then the estimate
$$\sum\limits_{R\le|z_j|<2R}r_j > R^{1-\varepsilon},\, R>
R(\varepsilon), $$ holds.
  \end{cit}
 \begin{cit}Let numbers $\rho >0,\,\varepsilon>0$ and an entire function $f$ satisfy the inequality
 $$||z|^\rho-\log|f(z)||< C_5\log |z|,\,z\notin E,
 $$
  where $E\subset \bigcup_j\{z:|z-z_j|<r_j\},\,r_j\le |z_j|^{1-\rho/2+\varepsilon}.$
Then the estimate
$$\sum\limits_{R\le|z_j|<2R}r_j>R^{1+\rho/2-2C_5-\varepsilon},\,
R> R(\varepsilon), $$ holds.
  \end{cit}
  Let us consider the accuracy of Theorem 1. It follows from Theorem B that there exist an entire function
   $f$ and an exceptional set $E,$ satisfying  (8) with $\lambda(R)=R$, $\tau=1$ and $m_2(E\bigcap A(R,2R])=o(R^{2\rho-2C_5/a_1-2a_0/a_1}),\,R\to \infty$.  We draw the conclusion that the planar Lebesgue measure of the exceptional set in the annulus $A(R,2R]$ is a power function of $R$ with the exponent $\asymp
  C_5$. The dependance of the exponent on $\rho$ is not clear.

       \textbf{Proof of Theorem 1.}
 We begin with the idea of the proof. At first, we  prove that any disk of the form
 $D(a,\,|a|^{1-\rho/2}),$ where
$f(a)=0$, contains a quite big portion of the exceptional set, namely, the estimate
 \begin{equation} m_2(E\cap
D(a,|a|^{1-\rho/2}))>C_7|a|^\chi
     \end{equation}  holds. The proof of estimate (10) is the key point of the proof of Theorem 1.
     We follow arguments from  ~\cite{H}, a new approach is that we use  the theorem by Edrei and Fuchs  on the integral over a small set~\cite{I}, ~\cite{K}. Next, it
is proved that every disk with somewhat greater radius has
the same property without the demand on the center of the disk to be a
zero of the function $f$. More exactly, we prove that for every
    $b,\,|b|>R_1,$ and every
$\varepsilon > 0$ the inequality \begin{equation} m_2(E\cap
D(b,|b|^{1-\rho/2+\varepsilon/2}))>C_8|b|^\chi
     \end{equation} holds. To finish the proof, we put
sufficiently many nonoverlapping  disks with enlarged radius into
the annulus
 $A(R,2R)$. Comparing the areas of the annulus and the disks, we obtain estimate
 (9).

We denote $r(a):=|a|^{1-\rho/2},\, v(z):=|z|^\rho$. Let  $h(z,v)$ and
$h(z,\log|f|)$ be the minimal harmonic majorants respectively $v$
and $\log|f|$ in the disk $D(a,t),\,t\in [(1-\varepsilon/4)r(a),\,r(a)]$.
Under the conditions of Theorem 1 we prove the estimate $(0<\delta<2\pi)$
\begin{equation}
|h(z,v)-h(z,\log|f|)|\le C_{9}\log \lambda (|a|)+C_{10}\delta
|a|^\rho \log |a| \log \left(\frac {2\pi e}{\delta }\right),z\in D(a,\,\varepsilon.
/4r(a)).
\end{equation}

 By the Poisson-Jensen formula for the disk $D(a,t)$ we have
 $$
|h(z,v)-h(z,\log|f|)|=$$ $$=\frac 1
{2\pi}\left|\int\limits_0^{2\pi}\left(|\,a+te^{i\varphi}|^\rho-
\log|f(a+te^{i\varphi})|\right)\Re
\frac{te^{i\varphi}+z-a}{te^{i\varphi}-z+a}\,d\varphi\right|\le
 $$
 $$
 \le \frac 1
{2\pi}\int\limits_0^{2\pi}\left||\,a+te^{i\varphi}|^\rho-
\log|f(a+te^{i\varphi})|\right|\frac
{r(a)+r(a)\varepsilon/4}{(1-\varepsilon/4)r(a)-r(a)\varepsilon/4}\,d\varphi\le
 $$\begin{equation}\le\frac {1+\varepsilon}
{2\pi}\int\limits_0^{2\pi}\left||\,a+te^{i\varphi}|^\rho-
\log|f(a+te^{i\varphi})|\right|\,d\varphi.
\end{equation}
 We denote by $E(t,a)$ the set
 $$ \left \{\varphi \in [0,2\pi]: ||\,a+te^{i\varphi}|^\rho-
\log|f(a+te^{i\varphi})||\ge C_5 \log \lambda
(|\,a+te^{i\varphi}|)\right.>$$\begin{equation}\left.>C_{5}(1-\varepsilon)\log
\lambda (|a|)\right\}\end{equation}
 (the last inequalities in definitions (14) and  next follow from the properties
 of the function $\lambda\in\Theta$ and the restrictions on $t$ for all sufficiently large values $|a|$). Its complement
 $$ [0,\,2\pi]\setminus E(t,a):=\left \{\varphi \in [0,2\pi]: ||\,a+te^{i\varphi}|^\rho-
\log|f(a+te^{i\varphi})||<\right.$$ $$\left.< C_5 \log \lambda
(|\,a+te^{i\varphi}|)<C_5(1+\varepsilon)\lambda(|a|)\right \}.$$
now we continue estimate (13):
 $$
 |h(z,v)-h(z,\log|f|)|\le\frac {1+\varepsilon}
{2\pi}\int\limits_{E(t,a)}\left||\,a+te^{i\varphi}|^\rho-
\log|f(a+te^{i\varphi})|\right|\,d\varphi+
$$
 $$
 +\frac {1+\varepsilon}
{2\pi}\int\limits_{[0,2\pi]\setminus
E(t,a)}\left||\,a+te^{i\varphi}|^\rho-
\log|f(a+te^{i\varphi})|\right|\,d\varphi\le
 $$
 $$
 \le\frac {1+\varepsilon}
{2\pi}\int\limits_{E(t,a)}\left(|\,a+te^{i\varphi}|^\rho+
\log^+|f(a+te^{i\varphi})|+\log^-|f(a+te^{i\varphi})|\right)\,d\varphi+
 $$
 \begin{equation}+ ({1+\varepsilon})^2C_{5}\log \lambda (|a|).\end{equation}

 We use the theorem by Edrei and Fuchs ~\cite{I},~\cite[С. 58]{K}, which we cite for the reader's convinience.\begin{cit} Let $f$ be a meromorphic function
 , $k$
  and $\delta$ real numbers, $k>1,\,0<\delta<2\pi,\,r>1$.
   For any measurable set $E_r\subset [0,2\pi]$, such that
     $m_1(E_r)\le\delta$, the relation
    \begin{equation}\int    \limits_{E_r}
    \log^+|f(re^{i\varphi})|\,d\varphi\le \frac{6k}{k-1}\,
    \delta\log\frac{2 \pi e}{\delta}\,T(kr,f)\end{equation}holds.\end{cit}
     We note that analysis of the proof of Theorem E shows that its
      $\delta$ - subharmonic version is valid. The assumption
     $r>1$ is of technical character: without it the term
      $$
    \delta\log\frac{2\sqrt {k}\pi e}{\delta}\,n(0,f)|\log r|/\log \sqrt{k}.
      $$
 should be added to the right-hand side (16)   (see the proof of Lemma 7.1 in
     ~\cite[P. 55]{K}). For the completeness of our exposition we perform a proof of above mentioned
            modification of Theorem E. We start from the inequalities  $(R'>R>0)$
$$N(R')\ge
\int\limits_R^{R'}\frac{n(t,f)-n(0,f)}{t}\,dt+n(0,f)\log R'\ge
(n(R,f)-n(0,f))\log \frac{R'}{R}+n(0,f)\log R,
$$
  from which the estimate
      $$ n(R,f)\le\frac{N(R',f)}{\log \frac{R'}{R}}-\frac{n(0,f)\log R}{\log
      \frac{R'}{R}}
      $$
      follows. In the proof of Lemma 7.1 it is supposed
   $R>1$, and because of this the negative term $-n(0,f)\log R$ is omitted, here
       we account it. After that in the proof of theorem E in
      ~\cite{K} the term
      $$n(\sqrt k r,f)\delta\log\frac{2\sqrt k e\pi}{\delta}
      $$
is obtained.
      To estimate $n(\sqrt k r,f)$ we apply the previous inequality with $R'=kr, R=\sqrt k r$,
      then we obtain
      $$n(\sqrt k r,f)\delta\log\frac{2\sqrt k e\pi}{\delta}\le \delta \log\frac{2\sqrt k e\pi}{\delta}
      \left(\frac{N(kr)}{\log \sqrt k}-\frac{n(0,f)\log (\sqrt k r)}{\log \sqrt
      k}\right)\le
      $$
      $$\le \delta \log\frac{2\sqrt k e\pi}{\delta}\left(\frac{N(kr)}{\log \sqrt k}-\frac{n(0,f)\log
       r}{\log \sqrt
      k}\right)\le \delta \log\frac{2\sqrt k
      e\pi}{\delta}\left(\frac{N(kr)}{\log \sqrt k}+\frac{n(0,f)|\log
       r|}{\log \sqrt
      k}\right),$$
      then we follow the proof in ~\cite{K}.

We continue estimate (15). The integral \begin{equation}
\int\limits_{E(t,a)}|a+te^{i\varphi}|^\rho\,d\varphi\le
(1+\varepsilon)|a|^\rho\delta
\end{equation} if $|a|$ is sufficiently large. To estimate the integral
$$\int\limits_{E(t,a)}\left(
\log^+|f(a+te^{i\varphi})|+\log^-|f(a+te^{i\varphi})|\right)\,d\varphi
$$ we apply precised Theorem E, putting $k=2$, and take into account the relation
 $T(r,f)=T(r,1/f)+O(1)$ and (7). We obtain the estimate
 $$\int\limits_{E(t,a)}\left(
\log^+|f(a+te^{i\varphi})|+\log^-|f(a+te^{i\varphi})|\right)\,d\varphi=$$
$$=\int\limits_{E(t,a)}\left( \log^+|f(a+te^{i\varphi})|+\log^+\left|\frac
 1
 {f(a+te^{i\varphi})}\right|\right)\,d\varphi \le $$
$$\le12\delta\log\frac{2 \pi e}{\delta}(2T(2t,f(z+a))+O(1))+
\delta\log\frac{2\sqrt {2}\pi e}{\delta}\,n(0,f(z+a))|\log t|/\log
\sqrt{2}< $$
\begin{equation}< 24 \delta\log\frac{2\pi e}{\delta}3^\rho C_4 |a|^\rho+6\delta\log\frac{2\sqrt {2}\pi
e}{\delta}\,C_4 2^\rho |a|^\rho |1-\rho/2| \log |a|.\end{equation}

Combining (15), (17), and (18), we have
\begin{equation}
|h(z,v)-h(z,\log |f|)|\le C_{11}\delta\log\frac{2\pi
e}{\delta}|a|^\rho \log |a|+C_5 (\tau+\varepsilon)\log |a|
\end{equation} if $|a|$ is sufficiently large
and $z\in D(a,\varepsilon r(a)/4).$

The next step is ti find an upper bound of the difference
 $\log |f(z)|-h(z,\log|f|)$ for $z\in D(a,\varepsilon r(a)/4)
\backslash E.$ It succeeds to prove such an estimate only in indirect way.
Using the standard tools of calculus, we prove
(see~\cite {H}) that for $z\in A:=A(R-r(R),R+r(R)]$ the inequality
\begin{equation}|v(z)-h(z,v,A)|\le C_{12}\end{equation} holds,
where $$h(z,v,A):=(R+r(R))^\rho\frac{\log
|z|-\log(R-r(R))}{\log(R+r(R))-\log(R-r(R))}+
$$
$$+(R-r(R))^\rho\frac{-\log
|z|+\log(R+r(R))}{\log(R+r(R))-\log(R-r(R))}
$$
is the minimal harmonic majorant of the function
$v(z):=|z|^\rho$ in the annulus $A$. From (20) and the definition of the minimal harmonic majorant it implies
\begin{equation}|v(z)-h(z,v,D(a,t))|\le C_{12}.\end{equation}
Applying (19), (21), and (8), we obtain the estimate $$ |\log
|f(z)|-h(z,\log|f|)|\le |h(z,v)-h(z,\log|f|)|+|v(z)-h(z,v)|+$$
$$+|v(z)-\log|f(z)||\le C_{11}\delta\log\frac{2 \pi
e}{\delta}|a|^\rho \log |a|+C_5 (\tau+\varepsilon)\log |a|+C_{12}+
$$
\begin{equation}+C_5 (\tau+\varepsilon)\log |a|\end{equation} if  $z\in D(a,\varepsilon r(a)/4)\setminus E.$

Now we prove the estimate from below of the difference $|\log
|f(z)|-h(z,\log|f|)|$. By the Poisson-Jensen formula
\begin{equation}\log
|f(z)|=h(z,\log |f|,D(a,\,t))-\sum\limits_{a_n\in D(a,\,t)}g(z,a_n),
\end{equation} where $g(z,a_n)$ is the Green function of the disk $D(a,\,t)$ with pole in zero $a_n$
of the function $f$. Using the known properties of the Green function, from (23)
we obtain
\begin{equation}|\log
|f(z)|-h(z,\log |f|,D(a,t))|\ge g(z,a)=\log\frac{t}{|z-a|}.
\end{equation}

We face the alternative: for every $t\in [(1-\varepsilon)r(a),r(a)]$
the measure  $m_1(E(t,a))\ge \delta$ or there exists $t\in
[(1-\varepsilon)r(a),r(a)]$, for which the measure $m_1(E(t,a))<
\delta$, where $\delta = \varepsilon (|a|^\rho \log |a|)^{-1}$. In the first case
the planar Lebesgue measure
\begin{equation}m_2(E\cap D(a,r(a)))\ge \varepsilon r(a)\cdot\varepsilon r(a)
(|a|^\rho \log |a|)^{-1}=\varepsilon^2\frac{|a|^{2-2\rho}}{\log
|a|}.\end{equation}
 In the second case, as it follows from (24), for $z\in D(a,(1-\varepsilon)r(a)/|a|^\kappa),\,
 \kappa=2C_5(\tau+\varepsilon)+2C_{11}\varepsilon,$ the estimate
   \begin{equation}|\log |f(z)|-h(z,\log|f|)| \ge \kappa \log |a|\end{equation} takes place, and from
    (22) we obtain that
    \begin{equation}|\log |f(z)|-h(z,\log|f|)| \le C_{11}\varepsilon\log |a|+
    2C_5(\tau+\varepsilon)\log |a|\end{equation}
   if $z\in D(a,\varepsilon r(a)/4)\setminus E.$ Comparing (26) and (27),
   we conclude, that the disk $D(a,(1-\varepsilon)r(a)/|a|^\kappa)\subset
   E.$
   Since its area $\pi (1-\varepsilon)^2 |a|^{2-\rho-2\kappa}$,
   and the planar Lebesgue measure of the portion  of the exceptional set $E\cap D(a,r(a))$
   does not exceed $\varepsilon^2|a|^{2-2\rho}/\log^2
|a|$ in the first case, then in any case (it is clear that $ \mathrm{const} \cdot \varepsilon $ can be replaced by $\varepsilon$)
\begin{equation}m_2(D(a,r(a))\cap E)\ge C_{13}|a|^\chi,\,\chi=\min
(2-2\rho-\varepsilon,2-\rho-4C_5\tau-\varepsilon).\end{equation}

  We put $r_1(b):=r(b)|b|^{\varepsilon/2}.$ For an arbitrary disk of the form
$D(b,r_1(b))$ with sufficiently large $|b|$ we prove the estimate
\begin{equation}m_2(D(b,r_1(b))\cap E)\ge C_{14}|b|^\chi.\end{equation}
Without lost of generality we can suppose $b\notin E$. In the opposite case
 $D(b,r(b))\subset E$ and then
(29) holds, or there exists $c\in D(b,r(b))\setminus E$. It is important that the disk
$D(c,\frac{3}{4}r_1(c))\subset D(b,r_1(b))$, this is applied under.
Suppose the disk $D(c,\frac{3}{4}r_1(c))$ does not contain zeros of the entire function
$f$, then $n(c,t,\log|f|)=0$ for $t\in [0,\frac{3}{4}r_1(c)]$.

 The Poisson-Jensen formula for the difference $|z|^\rho-\log|f(z)|$
in the disk $D(c,t/2)$ has the form
$$
-|c|^\rho+\log|f(c)|+\frac{1}{2\pi}\int\limits_{0}^{2\pi}\left(\left|c+\frac{1}{2}
te^{i\varphi}\right|^\rho-\log\left|f\left(c+\frac{1}{2}te^{i\varphi}\right)\right|\right)\,d\varphi=
$$
\begin{equation}=\int\limits_0^{t/2}\frac{n(c,s,|z|^\rho)}{s}\,ds.
\end{equation}

Estimating the integral on the left-hand side (30) like the integral in  (13), we obtain
$$\left|\frac{1}{2\pi}\int\limits_{0}^{2\pi}\left(\left|c+\frac{1}{2}
te^{i\varphi}\right|^\rho-\log\left|f(c+\frac{1}{2}te^{i\varphi})\right|\right)\,d\varphi\right|\le
$$
\begin{equation}\le C_{11}\delta\log\frac{2\pi
e}{\delta}T(t,f(w+c))+C_5(\tau+\varepsilon)\log
|c|+C_{15}\delta|c|^\rho,
\end{equation}
 where $m_1(E(t/2,c))\le\delta.$
Next, for $t\in [0,\frac{3}{4}r_1(c)]$ the estimate
\begin{equation}\frac{1}{2\pi}\int\limits_{0}^{2\pi}\log^+\left|f\left(c+\frac{1}{2}
te^{i\varphi}\right)\right|\,d\varphi\le C_4(|c|+t/2)^\rho \le
C_42^\rho |c|^\rho\end{equation} takes place if $|c|$ is sufficiently large. Here we make use of (7).
From (32) and the definition of the Navanlinna characteristic of a meromorphic function we obtain the estimate
\begin{equation}T(t,f(w+c))\le 2^{\rho }C_4
|c|^\rho.\end{equation} We put $\delta:=|c|^{-\rho }$. Again we face the alternative:
for every $t\in
[\frac{1}{2}r_1(c),\frac{3}{4}r_1(c)]$ the measure $m_1(E(t/2,c))\ge
\delta$, or there exists  $t\in
[\frac{1}{2}r_1(c),\frac{3}{4}r_1(c)]$ for which
$m_1(E(t/2,c))< \delta$. In the first case $$m_2(E\cap
D(b,r_1(b)))\ge m_2\left(E\cap
D\left(c,\frac{3}{4}r_1(c)\right)\right)\ge
|c|^{-\rho}\frac{1}{4}r_1(c)\frac{1}{4}r_1(c)\asymp$$ $$\asymp
|c|^{2-2\rho+\varepsilon}\asymp |b|^{2-2\rho+\varepsilon}. $$ In the second case
 (31) and (33) imply the inequality
\begin{equation}\left|\frac{1}{2\pi}\int\limits_{0}^{2\pi}\left|c+\frac{1}{2}
te^{i\varphi}\right|^\rho-\log^+\left|f\left(c+\frac{1}{2}te^{i\varphi}\right)\right|\,d\varphi\right|\le
C_{16}\log|c|+C_5(\tau+\varepsilon)\log|c|.
\end{equation} Since $c\notin E$, therefore\begin{equation}\left||c|^\rho-\log|f(c)|\right|\le C_5(\tau+\varepsilon)\log|c|.\end{equation}

Combining  (34) and (35), we conclude that the riht-hand side (30) does not exceed
 $C_{17}\log|c|.$

On the other hand, for the function  $v(z)=|z|^\rho$ its Riesz measure
$d\mu_v(z)=\frac{1}{2\pi}\Delta v\asymp |z|^{\rho-2}\,dm_2(z)$ and because of this
  $n(c,s,v)\asymp |c|^{\rho-2}s^2,$ if $s\le \frac{3}{4}
r_1(c),$ and the right-hand side of (30), i.e.
$\int_0^{t/2}\frac{n(c,s,v)}{s}\,ds\asymp
|c|^{\rho-2}r_1(c)^{2}\asymp |c|^{\varepsilon}(t\ge
\frac{1}{2}r_1(c))$, we obtain the contradiction with the previous estimate .
We draw a conclusion that there exists a zero $a$ of the entire function $f$ such that
$a\in D(c,\frac{3}{4}r_1(c))$. If $|b|$ is a sufficiently large number,
then $D(a,r(a))\subset D(b,r_1(b)).$ In any case the measure
$$ m_2(E\cap D(b,r_1(b)))\ge C_{8}|b|^\chi. $$ To finish the proof of Theorem 1
we put in the annulus $A[R,2R)$
nonoverlapping disks $D(b,r_1(b)))$ at a rate $\asymp
\frac{R^2}{R^{2-\rho+\varepsilon}}=R^{\rho-\varepsilon}$. Their union
contains such a portion of the exceptional set $E$,
that $$m_2(E\cap A[R,2R)) \ge C_6(\varepsilon)
R^{\chi+\rho-\varepsilon}.$$

Theorem 1 is proved.

I would like to thank all the referees of this work.

\end{document}